\documentclass[11pt]{amsart}
\usepackage{amsmath}
\usepackage{amscd}
\usepackage{amssymb}
\usepackage{graphicx}
\usepackage{graphics}
\usepackage{latexsym}
\usepackage{url}

\input xy
\xyoption{all}

\theoremstyle{definition}

\newcommand{\zr}{{\mathbb R}}

\begin{document}

\title{Theoretical Geometry, Critical Theory, and Concept Spaces in IR}
\author{Laura Sjoberg}
\email{sjoberg@ufl.edu}
\address{Department of Political Science, University of Florida, Gainesville, FL 32611}
\author{Kevin P.~Knudson}
\email{kknudson@ufl.edu }
\address{Department  of Mathematics, University of Florida, P.O.~Box 118105, Gainesville, FL 32611-8105}

\date{April 5, 2015}

\begin{abstract}  We use the theory of persistent homology to analyze a data set arising from the study of various aspects of democracy.  Our results show that most ``mature" democracies look more or less the same, in the sense that they form a single connected component in the data set, while more authoritarian countries cluster into groups depending on various factors.  For example, we find several distinct $2$-dimensional homology classes in the set, uncovering connections among the countries representing the vertices in the representative cycles.
\end{abstract}

\maketitle

As was discussed in the introduction to this book, the ‘quantitative’ methods traditionally used in the social sciences represent a limited subset of available methods in mathematics, statistics, and computational analysis, and the positivist ends for which they are usually deployed in the social science community represent a limited subset of the purposes for which they are intended and deployed in the philosophy of mathematics.

If it were to be oversimplified for explanatory purposes, math is the study of patterns. Discrete (if such things exist) ones are arithmetic.\footnote{E.g., Aristotle. See discussion in Jonathan Lear, “Aristotle’s Philosophy of Mathematics,” The Philosophical Review 91(2) (1982):161-192.}   Continuous ones are geometry.\footnote{E.g., Euclid. See discussion in J. C. Fisher, “Geometry According to Euclid,” The American Mathematical Monthly 86(4) (1979):260-270.}    Immeasurable ones are symbolic logic.\footnote{E.g., Bertrand Russell, “Mathematical Logic Based on the Theory of Types,” 
American Journal of Mathematics 30(3) (1908):222-262.}   Patterns not dependent on the empirical existence of their component parts are formalist.\footnote{E.g., Carl B. Boyer, “Descartes and the Geometrization of Algebra,” The American Mathematical Monthly 66(5) (1959):390-393.}  In the philosophy-of-math sense, this chapter takes a formalist approach to mathematical symbolism.\footnote{With Alan Weir (in . 2011. “Formalism in the Philosophy of Mathematics,” Stanford Encyclopedia of Philosophy [2011] accessed 15 March 2014 at http://plato.stanford.edu/entries/formalism-mathematics/), such an approach does not see mathematics as “a body of propositions representing an abstract sector of reality,” but “much more akin to a game, bringing with it no … commitment to an ontology of objects or properties.” In this sense, formalists “do not imply that these formulas are true statements” but instead see that “such intelligibility as mathematics possesses derives from the syntactical or metamathematical rules governing those marks.” (Nelson Goodman and W. V. Quine,  “Steps toward a Constructive Nominalism,” Journal of Symbolic Logic 12 [1947]:97-122, 122, 111).}

For the purpose of this analysis, that means that we argue that mathematical work is axioms with rules of inference that make possible thought experiments and string manipulation games of almost infinite complexity. There is no ‘true meaning’ underlying mathematical symbols – instead equations, formalizations, and quantifications are representations from which we learn about relationships  -homologies, homeomorphisms, maps, dimensionality, commutativity, factorization.\footnote{   In order to understand that assertion, it is important to distinguish, in discipline and in substance, between statistics and mathematics. Statistics is representation of quantity and correlation; mathematics is the spatial or equative representation of relationships, particularly as they become increasingly (and even impossibly) complex.  We contend that critical theory wants that, and math can do it. Math can do spatial concept maps of multidimensional, overlayered, folded, twisted manifolds of the stuff of politics. It can represent the mess. It can even suggest tweaks, trades, changes, foldovers – alterations to the map. Imagine math like the Matrix – but used to understand transparently rather than   manipulate. It is possible to look at this theoretical discussion as abstract, and impossible to materialize. That’s an admitted weakness.} In this view, the rules, laws, and procedures of mathematics are socially constructed – interested not only in quantity but also in structure, space, change, stochasticity, relationality, and formalization for its own sake. There are, for sure, many practitioners of math who think of it as science, as discovery, and as progress; but there are also many practitioners of math who see it as art, as creation, as signification, and as representation.\footnote{e.g., Nathalie Sinclair  “The Role of the Aesthetic in Mathematical Inquiry,” Mathematical Learning and Thinking 6(3) (2004):261-284; Edward A. Silver and Wendy Metzger, “Aesthetic Influences of Expert Mathematical Problem-Solving,” in Douglas B. McLeod and Verna M. Adams, eds. Affect and Mathematical Problem-Solving (New York: Springer, 1989), 59-74; Henri Poincare, “Mathematical Creation,” The Monist Volume XX (1910); Gontran Ervynckm “Mathematical Creativity,” Advanced Mathematical Thinking 11 (1991): 42-53; Loren Graham and Jean-Michael Kantor, Naming Infinity: A True Story of Religious Mysticism and Mathematical Creativity (Cambridge, MA: Harvard University Press, 2009); Gian-Carlo Rota, “The Phenomenology of Mathematical Beauty,” Synthese 111(2) (1997): 171-182; David Wells, “Which is the Most Beautiful,” The Mathematical Intelligencer 10(4) (1988): 30-31. }

The overwhelming majority of quantitative methods that are used in International Relations (IR) are statistical in nature. When the tools of theoretical mathematics are used, they are most often deployed in predictive, descriptive, or heuristic uses of game theory. That is, the tools of mathematics are used to try to gain leverage on the causal empirical realities that neopositivist IR scholars see as the important substance of global politics to know and understand. In this chapter, we argue that, not only can the tools of theoretical mathematics be utilized for post-positivist ends, but that is where many of those tools would be most at home in IR. 

Particularly, in this chapter, we look to make an initial case for the argument that the tools of computational topology can be effectively utilized to explore questions of constitution, textuality, and performativity for critical IR.  To that end, the first section lays out an argument about the possible utility of thinking geometrically about concept formation and reification for post-structuralist IR, and possible ways to do that work. The second section introduces the concept of democracy in IR, and argues that it might be possible to gain leverage on the dimensions of the concept using computational topology to evaluate existing data. The third section shows the method in action, and sketches out some of the possible ramifications for studying democracy from a critical perspective. The concluding section makes a case for the value-added both for political methodology and critical theory of methodological explorations like this. 

\section*{Critical IR and Geometry Topology}

In Mark Hoffman’s words, post-structuralist IR asks “how” questions, rather than “what” or “why” questions.\footnote{Mark Hoffman, “Critical Theory and the Interparadigm Debate,” Millennium: 
Journal of International Studies 16 (1987):231-249.} 
  Those “how” questions include: “how are structures and practices replicated? How is meaning fixed, questioned, reinterpreted, and refixed?”\footnote{Ibid.}   To answer these “how” questions, post-structuralist IR rejects “the modernist belief in our ability to rationally perceive and theorize the world” in favor of “dis-belief in unproblematic notions of modernity, enlightenment, truth, science, and reason.”\footnote{Lene Hansen,  “A Case for Seduction? Evaluating the Poststructuralist Conceptialization of Security,” Cooperation and Conflict – Nordic Journal of International Studies 34(2) (1997): XXXX.} This move leads post-structuralists, methodologically, to a sort of scholarship which “does not look for a continuous history, but for discontinuity and forgotten meanings; it does not look for an origin, indeed, it is assumed one cannot be found; and it does not, finally, focus on the ‘object of geneaology’ itself, but on the conditions, discourses, and interpretations surrounding it.”\footnote{Ibid., p.372.}  This scholarship pays particular attention to discourses in global politics, where (actual and conceptual) “boundaries are constantly being redrawn and transgressed.”\footnote{Hoffman, “Critical Theory and the Interparadigm Debate.”}

We argue that this {\em methodological} approach can be matched fairly easily with the capacities of the {\em methods} of geometric and computational topology. If critical IR rejects the objective and the linear, and with them their dichotomized frame of reference, topological analysis can accommodate significantly more complexity. If critical IR tries to unmask and deconstruct hidden meanings, there could be some power in representing those meanings geometrically. If critical IR embraces undecidability, is at home with liminality, is wary of metanarratives, and is attached to the “how” questions discussed above, the instability, creativity, and formalism of theoretical geometry may be a good fit.    

In fact, this is not the first time a critical theorist has suggested that there might be some benefit to thinking geometrically about the concepts being explored and critiqued. Deleuze and Guattari see concepts as rhizomes, biological entities endowed with unique properties.\footnote{Gilles Deleuze and Felix Guattari, A thousand plateaus: capitalism and schizophrenia,  translated by Brian Massumi.  (Minneapolis: University of Minnesota Press, 1987).}  They see concepts as spatially representable, where the representation contains “principles of connection and heterogeneity: any point of a rhizome must be connected to any other.”\footnote{Ibid., p. 7. As the authors continue on to explain: “A rhizome ceaselessly establishes connections between semiotic chains, organizations of power, and circumstances relative to the arts, sciences, and social struggles.  A semiotic chain is like a tuber agglomerating very diverse acts, not only linguistic, but also perceptive, mimetic, gestural, and cognitive: there is no language in itself, nor are there any linguistic universals, only a throng of dialects, patois, slangs, and specialized languages.”}   Deleuze and Guattari list the possible benefits of spatial representation of concepts, including the ability to represent complex multiplicity, the potential to free a concept from foundationalism, and the ability to show both breadth and depth.\footnote{Ibid., p.8.}   In this view, geometric interpretations move away from the insidious understanding of the world in terms of dualisms, dichotomies, and lines, to understand conceptual relations in terms of space and shapes. The {\em ontology} of concepts is thus, in their view, appropriately geometric – a multiplicity “defined not by its elements, nor by a center of unification and comprehension” and instead measured by its dimensionality and its heterogeneity.\footnote{Ibid., p.267.}  The conceptual multiplicity, “is already composed of heterogeneous terms in symbiosis , and … is continually transforming itself” such that it is possible to follow, and map, not only the relationships between ideas but how they change over time.\footnote{Ibid.}  In fact, the authors claim that there are further benefits to geometric interpretations of understanding concepts which are unavailable in other frames of reference.  They outline the unique contribution of geometric models to the understanding of contingent structure: ``Principle of cartography and decalcomania: a rhizome is not amenable to any structural or generative model.  It is a stranger to any idea of genetic axis or deep structure.  A genetic axis is like an objective pivotal unity upon which successive stages are organized; deep structure is more like a base sequence that can be broken down into immediate constituents, while the unity of the product passes into another, transformational and subjective, dimension."  (Deleuze and Guattari 1987, 12)

The word that Deleuze and Guattari use for `multiplicities' can also be translated to the topological term `manifold.'\footnote{If we thought about their multiplicities as manifolds, there are a virtually unlimited number of things one could come to know, in geometric terms, about (and with) our “object of study,” abstractly speaking. Among those unlimited things we could learn are properties of groups (homological, cohomological, and homeomorphic), complex directionality (maps, morphisms, isomorphisms, and orientability), dimensionality (codimensionality, structure, embeddedness), partiality (differentiation, commutativity, simultaneity), and shifting representation (factorization, ideal classes, reciprocity). Each of these functions allows for a different, creative, and potentially critical representation of global politics’ concepts, events, groupings, and relationships.}  This is how we propose looking at concepts: as manifolds. With such a dimensional understanding of concept-formation, it is possible to deal with complex interactions of like entities, and interactions of unlike entities.  Critical theorists have emphasized the importance of such complexity in representation a number of times, speaking about it in terms compatible with mathematical methods if not mathematically. For example, Michel Foucault’s declaration that ``practicing criticism is a matter of making facile gestures difficult” both reflects and is reflected in many critical theorists’ projects of revealing the complexity in (apparently simple) concepts deployed both in global politics and IR scholarship.\footnote{Michel Foucault, Politics, Philosophy, Culture (New York: Routledge, 1988), p.154.}  David Campbell's reading of the state in Writing Security is a good example of this: Campbell makes the argument that the notion of the state appears to be both simple and {\em a priori}, it is really danger built over other danger--where ``the constant articulation of danger through foreign policy is thus not a threat to a state’s identity or existence: it is its condition of possibility.”\footnote{David Campbell, Writing Security: United States Foreign Policy and the Politics of Identity (Minneapolis: University of Minnesota Press, 1992), p.10.}  This leads to a shift in the concept of danger as well, where ``danger is not an objective condition” but ``an effect of interpretation.”\footnote{Ibid., p.1, 2.}  Critical thinking about how-possible questions reveals a complexity to the concept of the state which is often overlooked in traditional analyses, sending a wave of added complexity through other concepts as well. This work {\em seeking complexity} serves one of the major underlying functions of critical theorizing: finding invisible injustices in (modernist, linear, structuralist) givens in the operation and analysis of global politics. 

In a geometric sense, this complexity could be thought about as multidimensional mapping. In theoretical geometry, the process of mapping conceptual spaces is ``not primarily empirical”\footnote{Peter Gardenfors, “Mental Representation, Conceptual Spaces, and Metaphors,” Synthese 106(1) (1996): 21-47.} but for the purpose of representing and reading the relationships between information, including identification, similarity, differentiation, and distance.\footnote{Janet Aisbett and Greg Gibbon, “A General Formulation of Conceptual Spaces as Meso-Level Representation,” Artificial Intelligence, 133 (2001): 189-232, p.190.} The reason for defining topological spaces in math, the essence of the definition, is that there is no absolute scale for describing the distance or relation between certain points, yet it makes sense to say that an (infinite) sequence of points approaches some other (but again, no way to describe `how quickly’ or `from what direction’ one might be approaching).  This seemingly weak relationship, which is defined purely `locally’, i.e., in a small locale around each point, is often surprisingly powerful: using only the relationship of `approaching parts’, one can distinguish between, say, a balloon, a sheet of paper, a circle, and a dot.

To each delineated concept, one should distinguish and associate a topological space, in a (necessarily) non-explicit yet definite manner.  Whenever one has a relationship between concepts (here we think of the primary relationship as being that of constitution, but not restrictively, we `specify' a function (or inclusion, or relation) between the topological spaces associated to the concepts). In these terms, ``a conceptual space is in essence a multidimensional space in which the dimensions represent qualities or features of that which is being represented.”\footnote{Ibid., p.192.}   Such an approach can be leveraged for thinking about conceptual components, dimensionality, and structure.\footnote{Gardenfors, “Mental Representation,” p.33.}   In these terms, dimensions can be thought of as properties or qualities, each with their own (often-multidimensional) properties or qualities.\footnote{Aisbett and Gibbon, “A General Formulation of Conceptual Spaces,” p.192.}

A key goal of the modeling of conceptual space being representation means that a key (mathematical and theoretical) goal of concept space mapping is  ``associationism, where associations between different kinds of information elements carry the main burden of representation."\footnote{Peter Gardenfors, “Conceptual Spaces as a Framework for Knowledge Representation,” Mind and Matter 2(2) (2004): 9-27, p.9.}  To this end, ``objects in conceptual space are represented by points, in each domain, that characterize their dimensional values."\footnote{John T. Rickard, “A Conceptual Geometry for Conceptual Space,” Fuzzy Optimum Decision Making 5 (2006): 311-329 (original emphasis removed).} These dimensional values can be arranged in relation to each other, as Gardenfors explains that ``distances represent degrees of similarity between objects represented in space" and therefore conceptual spaces are ``suitable for representing different kinds of similarity relation."\footnote{Peter Gardenfors, “Conceptual Learning: A Geometrical Model,” Meeting of the Aristotelian Society, University of London, Senate House, 5 February 2001, p.169.}  These similarity relationships can be explored across ideas of a concept and across contexts, but also over time, since ``with the aid of a topological structure, we can speak about continuity, e.g., a {\em continuous change}" – a possibility which can be found {\em only} in treating concepts as topological structures and not in linguistic descriptions or set theoretic representations.\footnote{Thomas Mormann, “Natural Predicates and Topological Structures of Conceptual Spaces,” Synthese 95(2) (1993): 219-240, p.220.}  Such an approach is both complex and {\em anexact} – suiting it well for the contingent explorations of critical IR.    

\bigskip

\noindent {\em A Formalization of Concept Relationships}

\bigskip

The first step might be to gain information about the (actual, representational, or potential) relationship between a concept being examined and another concept that contributes something to the essence of how it is understood.  Assume a complex concept K composed of (but not necessarily limited to) component parts $v_0,v_1,\dots ,v_p$.  The concept can be explored as a simplicial homology, where an {\em abstract simplicial complex} $K$ is specified by the following data:
\begin{itemize}
\item[$\bullet$] A vertex set $V$;

\item[$\bullet$] A rule specifying when a $p$-simplex $\sigma=[v_0 v_1 \dots v_p]$ belongs to $K$; here the vertices $v_0, v_1,\dots ,v_p$ are distinct elements of $V$;

\item[$\bullet$] Each $p$-simplex $\sigma$ has $p+1$ faces which are the $(p-1)$-simplices obtained by deleting one of the vertices of $\sigma$.  The membership rule has the property that if $\sigma$ belongs to $K$, then all of its faces belong to $K$.
\end{itemize}

Given a simplicial complex $K$, we wish to define a collection of vector spaces which tell us the number of holes of various dimensions in $K$. Before defining these objects formally, let us consider some simple examples.  Let $T$ be a hollow tetrahedron; topologically, this is a sphere.  Observe that any simple closed loop on the surface of $T$ bounds a disc.  This means that there are no one-dimensional loops in $T$ that cannot be filled in by a two-dimensional surface in $T$.  There is a two-dimensional surface, though, namely $T$ itself, that cannot be filled in without leaving $T$.  The corresponding homology vector spaces would then be zero-dimensional in degree 1 (every loop bounds a disc), and one-dimensional in degree 2 (there is a two-dimensional surface in $T$ that is not filled in by a three-dimensional object).  By contrast, consider the surface $S$ of a donut.  This also has a three-dimensional void (the interior of $S$), but it also has some loops that do not bound discs on $S$---choose any meridian and any longitude.  The degree 1 homology space would then have dimension 2.

 The formal definition is as follows. Let $k$ be the field of 2 elements.  The $i$-th homology group, $H_i(K;k)$, will measure the number of $(i+1)$-dimensional voids, and it is constructed as follows.  Let $C_i(K;k)$ be the $k$-vector space with basis the set of $i$-simplices in $K$.  If $\sigma=[v_0 v_1 \cdots v_i]$ is such a simplex, we define $\partial\sigma$ to be the element of $C_{i-1}(K;k)$ given by the formula
$$\partial\sigma=\sum_{j=0}^i (-1)^j[v_0 v_1 \cdots \hat{v}_j\cdots v_i],$$ where $[v_0 v_1 \cdots \hat{v}_j\cdots v_i]$ is the $(i-1)$-simplex with vertices $\{v_0,\dots ,v_i\}-\{v_j\}$. Note that in the field of 2 elements we have $-1=1$, but we present the definition this way because it works over any field (e.g., the real numbers). We may extend this linearly to $C_i(K;k)$ to obtain a linear transformation $$\partial^i:C_i(K;k)\to C_{i-1}(K;k).$$ It is a straightforward exercise to show that $\partial^{i}\circ\partial^{i+1}=0$ and hence $\textrm{im}(\partial^{i+1})\subseteq\textrm{null}(\partial^{i})$ (null($\partial^i$) denotes the null space of the map $\partial^i$).  We then define the $i$-th homology group as
$$H_i(K;k) = \textrm{null}(\partial^i)/\textrm{im}(\partial^{i+1}).$$  Elements of $\textrm{null}(\partial^i)$ are called {\em cycles}; the set of all such is denoted by $Z_i$.  Elements of $\textrm{im}(\partial^{i+1})$ are called {\em boundaries}, denoted by $B_i$.  Homology measures how many cycles are inequivalent and essential in the sense that they do not bound an object of higher dimension.  Note that the group $H_0$ measures how many connected components the space $K$ has.

For computational purposes, the goal is to pay attention to the {\em Betti numbers}, $\beta_i$, defined as $\beta_i=\textrm{dim}_k H_i(K;k)$. Note that we have the simple equation
\begin{eqnarray*}
\beta_i & = & \dim_k Z_i - \dim_k B_{i+1} \\
        & = & \dim_k C_i - \textrm{rank}\,\partial^i - \textrm{rank}\,\partial^{i+1}.
\end{eqnarray*}
This therefore reduces the calculation of Betti numbers to computing ranks of matrices over the field $k$.

As a simple example, consider the tetrahedron $T$.  There are four vertices $v_0,v_1,v_2,v_3$; six edges $[v_0 v_1]$, $[v_0 v_2]$, $[v_0 v_3]$, $[v_1 v_2]$, $[v_1 v_3]$, $[v_2 v_3]$; and four faces $[v_0 v_1 v_2]$, $[v_0 v_1 v_3]$, $[v_0 v_2 v_3]$, $[v_1 v_2 v_3]$.  The groups $C_i(T;k)$ are then $$C_0(T;k) = k^4\qquad C_1(T;k)=k^6\qquad C_2(T;k)=k^4$$ and the maps $\partial^i$ are given by $\partial^0 = 0$, $\partial^i=0, i\ge 3$, and
$$\partial^1 = \left[\begin{array}{cccccc}
                   -1 & -1 & -1 & 0 & 0 & 0 \\
                   1 & 0 & 0 & -1 & -1 & 0  \\
                   0 & 1 & 0 & 1 & 0 & -1 \\
                   0 & 0 & 1 & 0 & 1 & 1
                   \end{array}\right]  $$
$$\partial^2 = \left[\begin{array}{cccc}
                      1 & 1 & 0 & 0 \\
                      -1 & 0 & 1 & 0 \\
                      0 & -1 & -1 & 0 \\
                      1 & 0 & 0 & 1 \\
                      0 & 1 & 0 & -1 \\
                      0 & 0 & 1 & 1
                      \end{array}\right].$$
An easy calculation shows that $\textrm{rank}\,\partial^1=3$, $\textrm{rank}\,\partial^2=3$, and therefore that
\begin{eqnarray*}
H_0(T;k) = k & \Rightarrow & \beta_0 = 1 \\
H_1(T;k) = 0 & \Rightarrow & \beta_1 = 0 \\
H_2(T;k) = k & \Rightarrow & \beta_2 = 1.
\end{eqnarray*}
A basis for $H_0(T;k)$ is $[v_0]$, and for $H_2(T;k)$, $[v_0 v_1 v_2] - [v_0 v_1 v_3] + [v_0 v_2 v_3] - [v_1 v_2 v_3]$.  There are three linearly independent $1$-cycles---the boundaries of the four triangles form a linearly dependent set of dimension $3$---but each is also a boundary, filled in by the interior of the triangle.  Geometrically, the fact that $\beta_0=1$ means that $T$ is connected; $\beta_1=0$ means that every loop in $T$ bounds a disc; $\beta_2=1$ means that $T$ contains a closed surface, namely $T$ itself, that is not filled in by a $3$-dimensional object.

Homology groups are topological invariants; that is, if spaces $X$ and $Y$ are homotopy equivalent (one may be deformed to the other), then $H_\bullet(X;k)\cong H_\bullet(Y;k)$.  They therefore provide a means to distinguish spaces, although it is possible for topologically distinct spaces to have the same homology groups.  

Tools for analyzing homology groups then can be used for similarity and difference analysis in these concept spaces. An increasingly popular technique for analyzing data sets topologically is the {\em persistent homology} of Edelsbrunner, Letscher, and Zomorodian.\footnote{Herbert Edelsbrunner; David Letscher; and Afra Zomordian, “Topological Persistence and Simplification,” Discrete \& Computational Geometry 28 (2002): 511-533.}  The idea is as follows.  Suppose we are given a finite nested sequence of finite simplicial complexes
$$K_{R_1}\subset K_{R_2}\subset \cdots \subset K_{R_p},$$ where the $R_i$ are real numbers $R_1<R_2<\cdots <R_p$.  For each homological degree $\ell\ge 0$, we then obtain a sequence of homology groups and induced linear transformations
$$H_\ell(K_{R_1})\to H_\ell(K_{R_2})\to\cdots\to H_\ell(K_{R_p}).$$  Since the complexes are finite, each $H_\ell(K_{R_i})$ is a finite-dimensional vector space.  Thus, there are only finitely many distinct homology classes.  A particular class $z$ may come into existence in $H_\ell(K_{R_s})$, and then one of two things happens.  Either $z$ maps to $0$ (i.e., the cycle representing $z$ gets filled in) in some $H_\ell(K_{R_t})$, $R_s<R_t$, or $z$ maps to a nontrivial element in $H_\ell(K_{R_p})$.  This yields a {\em barcode}, a collection of interval graphs lying above an axis parametrized by $R$.  An interval of the form $[R_s,R_t]$ corresponds to a class that appears at $R_s$ and dies at $R_t$.  Classes that live to $K_{R_p}$ are usually represented by the infinite interval $[R_s,\infty)$ to indicate that such classes are real features of the full complex $K_{R_p}$.

\begin{figure}
\centerline{\includegraphics[width=4in]{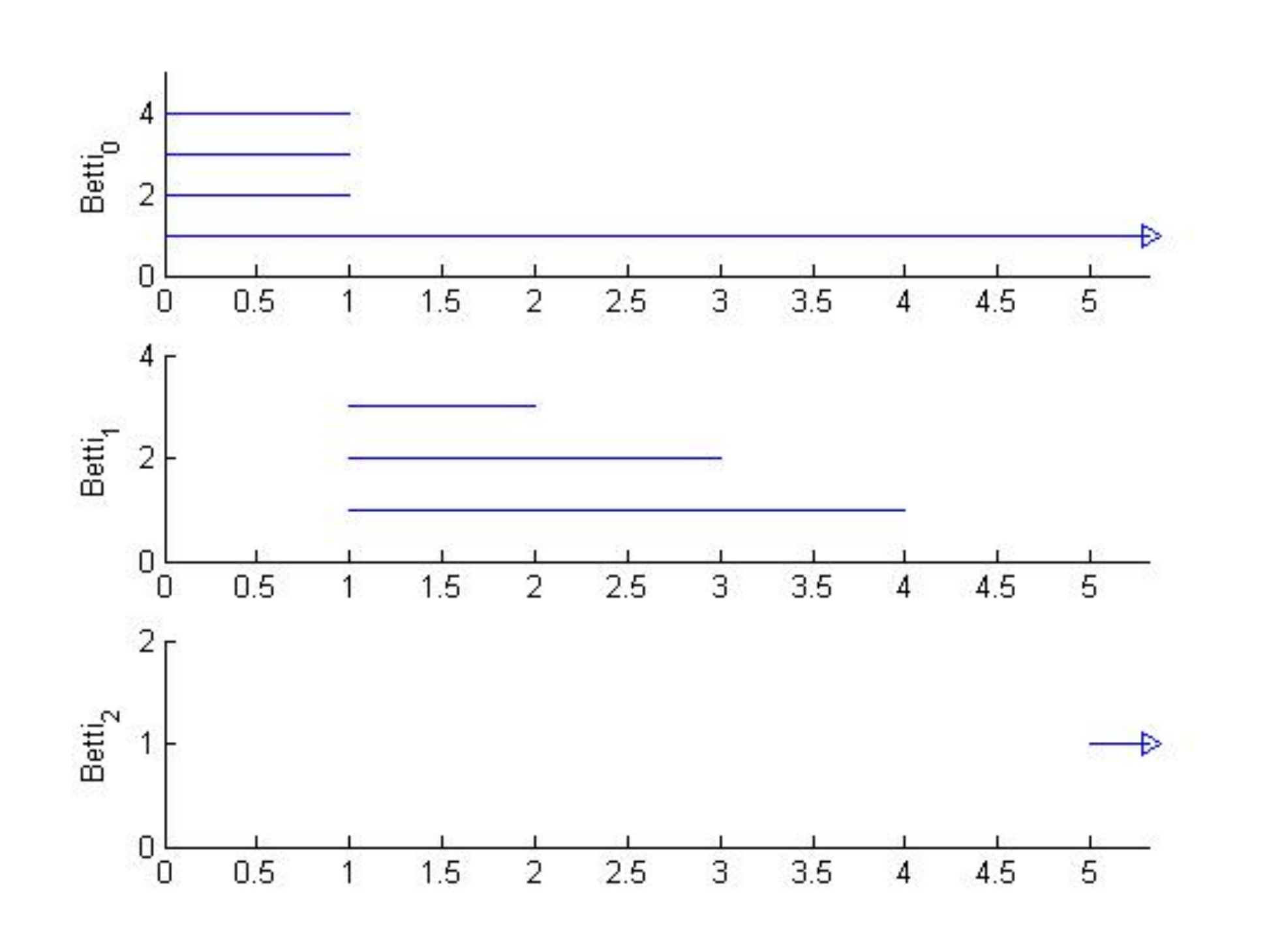}}
\caption{\label{tetrafilt} The barcodes for a filtration of the tetrahedron.}
\end{figure}

As an example, consider the tetrahedron $T$ with filtration $$T_0\subset T_1\subset T_2\subset T_3\subset T_4\subset T_5=T$$ defined by $T_0=\{v_0,v_1,v_2,v_3\}$, $T_1=T_0 \cup \{\textrm{all edges}\}$, $T_2=T_1\cup [v_0 v_1 v_2]$, $T_3=T_2\cup [v_0 v_1 v_3]$, $T_4=T_3\cup [v_0 v_2 v_3]$, and $T_5=T$.  The barcodes for this filtration are shown in Figure \ref{tetrafilt}.  Note that initially, there are $4$ components ($\beta_0=4$), which get connected in $T_1$, when $3$ independent $1$-cycles are born ($\beta_1=3$).  These three $1$-cycles die successively as triangles get added in $T_2$, $T_3$, and $T_4$.  The addition of the final triangle in $T_5$ creates a $2$-cycle ($\beta_2=1$).  We see that the classes that `live forever' yield the Betti numbers discussed above.

This approach can be used to sketch out relationships among concepts without data about how they are traditionally considered, or to analyze the constitution of concepts for which point cloud data is collected. For analyzing point cloud data, one needs a simplicial complex modeling the underlying space.  Since it is impossible to know {\em a priori} if a complex is ``correct", one builds a nested family of complexes approximating the data cloud, computes the persistent homology of the resulting filtration, and looks for homology classes that exist in long sections of the filtration.    A popular technique for this is the  {\em witness complex} construction of de Silva and Carlsson.\footnote{V. de Silva and G. Carlsson, “Topological Estimation Using Witness Complexes,” Symposium on Point-Based Graphics, ETH, Zurich, Switzerland, June 2-4, 2004.} 

Suppose there is a data set $Z$.  Choose a subset $L\subseteq Z$; we call the elements of $L$ {\em landmarks}.  Say $Z$ has $N$ elements and $L$ has $n$ elements and denote by $D$ the $n\times N$ matrix of distances between the elements of $L$ and the elements of $Z$.  Fix a real number $R>0$.  The {\em witness complex} $W(D,R)$ is defined as follows.

\begin{enumerate}
\item The vertex set of $W(D,R)$ is $\{1,2,\dots ,n\}$.
\item The edge $\sigma=[ab]$ belongs to $W(D,R)$ if and only if there exists a data point $1\le i\le N$ such that $$\max (D(a,i),D(b,i))\le R.$$  In this case the point $i$ is called a {\em witness} for $\sigma$.
\item The $p$-simplex $\sigma=[a_0a_1\dots a_p]$ belongs to $W(D,R)$ if and only if all its edges belong to $W(D)$; equivalently, there exists a witness $1\le i\le N$ such that $$\max (D(a_0,i),D(a_1,i),\dots ,D(a_p,i)) \le R.$$
\end{enumerate}

This definition may seem a bit opaque, but can be described geometrically as follows.  Imagine the data set $Z$ lies in some Euclidean space.  We then have the collection $L$ of landmarks and we build a simplicial complex with vertex set $L$ by looking for points in the data set near the landmarks.  Two landmarks are joined by an edge precisely when there is a data point in $Z$ within $R$ of both of them.  In fact, if we take $L=Z$, then this process yields the {\em Rips complex} on $Z$ (with parameter $R/2$) in which two data points are joined when their distance apart is at most $R/2$.  A good image to keep in mind for the Rips complex is that we place small balls of radius $r>0$ around each data point.  We join two data points if the balls meet.  That this is a good thing to do from a topological point of view is well-established.\footnote{See, e.g., James Munkres, Elements of Algebraic Topology (Menlo Park, CA: Benjamin/Cummings, 1984).}  

In the examples we analyze in this paper, we take $L=Z$ and therefore compute the Rips complex on the underlying data set.  For larger data sets, however, it is often impractical to build a simplicial complex on the entire set and one therefore turns to the witness complex construction on a much smaller subset of landmarks.  We refer the reader to de Silva and Carlsson (2004) for an analysis of the effectiveness of this process.

Note that the construction of $W(D,R)$ takes a part of its input the parameter $R$.  This allows us to construct a filtered collection of complexes by letting $R$ increase from $0$ to some large upper bound.  Note that if $R$ is sufficiently large, the complex $W(D,R)$ contains all possible simplices on the landmark set and is therefore a contractible $(n-1)$-simplex.  However, each simplex $\sigma$ comes into existence at some particular value $R_\sigma$, called its {\em time of appearance}.  Since there are only finitely many simplices possible on any landmark set, there is a discrete collection of parameter values $0<R_1<R_2<\cdots <R_r$ for which we get distinct witness complexes $$L=W(D,0)\subset W(D,R_1)\subset \cdots \subset W(D,R_r).$$  We may then compute the persistent homology of this filtration, searching for homology classes that persist for long intervals $[R_s,R_t]$.

In this way, analysis of persistent homologies can be used to show both similarities among concepts {\em and} similarities across components of concepts in particular cases to which those concepts and their components are applied. For critical IR, mapping conceptual spaces like this can provide a ``framework for representation" that demonstrates relationships among concepts without reaching necessary or essential conclusions about their genesis or origin (Gardenfors 1996). Analyzing the geometric complexity of concepts could lead to the ability to gain leverage on how understanding conceptual dimensionality could make current research better, more interesting, or deeper, either on its own terms or in terms of understanding complexity, hybridity, marginalization, social disadvantages, or other areas of global politics of interest to critical theory.  

\section*{Thinking Critically about Democracy}
Scholars in Comparative Politics and IR have long been interested in the nature and existence of democracy in the global political arena. Scholars of Comparative Politics investigate the structure and function of democratic institutions in states that they see as transitioning to democracy and in states they see as developing or mature democracies. Scholars of IR look to understand the ways in which states’ regime types may affect their foreign policy propensities, including but not limited to trade patterns, likelihood to be involved in conflict, and conflict opponents.
 
Scholars in Comparative Politics and International Relations make a number of common assumptions about democracy that permeate a significant amount of the research which is read, cited, and engaged in their respective fields. Each assumes that democracy is an extant and practiced form of government. Each assumes that such a form of government is {\em measurably} different than other forms of government, which may vary but share the label non-democratic. Though not all scholars make the assumption that democracy is a more desirable form of government than autocracy, oligarchy, theocracy, or other possible forms of government, many if not most of them do. 

That said, while scholars interested in either the internal or external politics of either democracies or democratizing states agree on the existence and distinguishability of democracy, many of them disagree on the components of democracy which are distinguishable, or on the specific places and times in which democracy exists. In other words, while the idea of democracy is common among scholars interested in comparative or international politics, many of them disagree on what makes a democracy, which countries are democracies, which components of democracy are measurable, and which measurable components of democracy are most central to the concept of democracy. 
	
Critical theorists have been interested in the concept of democracy in the international arena in a number of different ways for an extended period of time.  Of course, neither the breadth nor theoretical depth of critical analyses of democracy (even in IR) can be done justice in this small section of a chapter. That said, some of the common critiques (and resultant how-possible questions) can be explored briefly to give a sense of what the stake(s) in the concept of democracy may be for (different) critical IR researchers. Critical theorists have been concerned with the meaning of the concept of democracy, about the structure and normative value of the signification of the concept, and about the potential to revise and reappropriate the concept in search of a more just global political arena.  

Some critical theorists have been interested in the way that democracy as a concept signifies the success of `the West' and distinguishes that from an othered rest of the world.  As Amitav Acharya and Barry Buzan note, ``the contemporary equivalent of `good life' in international relations – democratic peace, interdependence and integration, and institutionalized orderliness, as well as the `normal relationships and calculable results' – is mostly found in the West, while the non-West remains the realm of survival."\footnote{Amitav Acharya and Barry Buzan, “Why Is There No Non-Western International Relations Theory? An Introduction,” International Relations of the Asia-Pacific 7(3) (2007): 287-312,  p.288, citing J. M. Goldgeier and M. McFaul, “A Tale of Two Worlds: Core and Periphery in the post-Cold War Era,” International Organization 46(2) (1992): 467-491.}  Accordingly Acharya and Buzan characterize democracy as a ``Western idea."\footnote{Ibid., p.295.}   As Fabrizio Eva notes, these origins hold in contemporary global politics, and ``the model for democracy nevertheless remains the Western version, which is tied in with capitalist economics. There is an acceptance, therefore, of the central ideology of liberal democracy."\footnote{Ibid., p.295.}  Many critical theorists’ concern with the Western-centric nature of deployed concepts of democracy in the international arena causes them (especially those interested in the question from a postcolonial perspective) to be critical of the use of the idea of democracy {\em writ large} in global politics or the analysis thereof.  For example, Inayatullah and Blaney argue that ``democracy is less a form of governance than a value that must be moderated, a set of practices to be disciplined by some prior claim to authority."\footnote{Naeem Inayatullah and David Blaney, International Relations and the Problem of Difference (New York and London: Routledge, 2004).}   Another reading of this analysis might suggest that democracy is not {\em a thing} out there to be analyzed, achieved, or deconstructed, but instead {\em a signifier} by which participants in the global political arena are organized hierarchically.  

Roxanne Doty sees this both in the policy arena and in the academic study of democracy.\footnote{Roxanne Lynn Doty, Imperial Encounters: The Politics of Representation in North-South Relations (Minneapolis: University of Minnesota Press, 1996), p.135.}  Doty explains that work in IR interested in democracy often presumes that ``some subjects were definers, delimiters, and boundary setters of important … and that others not capable themselves of making such definitions, would have things bestowed upon them and would be permitted to enjoy them only under the circumstances deemed suitable by the United States."\footnote{Ibid.}  This is why Tanji and Lawson argue that ``the `answer' to the question of what constitutes `true democracy' is implicit in the model of democracy assumed by the thesis … [which is] authoritatively assumed in advance, posited as an unassailable universal, and deployed as the foundation of the moral high ground in the global sphere."\footnote{Miyume Tanji and Stephanie Lawson, “’Democratic Peace’ and ‘Asian Democracy’: A Universalist-Particularist Tension,” Alternatives 22(1) (1997): 135-155, p.151. See discussion in Inayatullah and Blaney, The Problem of Difference.}   These readings suggest a critical stance not only towards current concepts of and wieldings of democracy, but also towards the use of the word and idea in general. This work in critical IR, then, is less concerned with reviving and rehabilitating the notion of democracy in global politics and more concerned with remedying the exclusions and silences produced by its current significations.

Still, even the most skeptical of critical theorists pay attention to the multiple significations of the term and idea of democracy in global politics. Doty draws on Laclau’s use of the signifier of democracy, explaining that it ``acquires particular meanings when it is associated with other signifiers," that is, that anticommunist democracy and antifascist democracy are signified differently, even if each is ``an attempt to constitute a hegemonic formation."\footnote{Doty, Imperial Encounters, pp. 8-9.}  Internally linked significations include democracy’s perceived opponents and its perceived components and benefits. For example, Doty suggests that American masculinity is a key component and tie-in to 19th century notions of democracy in global politics, where ``American manhood was also linked to democracy" and ``this link served to construct a distinctly American version of masculinity that was part and parcel of American exceptionalism."\footnote{Ibid., p.31.}   In other words, democracy and masculinity were co-constituted in a particular instantiation of democracy in global politics. More recently, Zalewski and Runyan suggest that the signification of democracy has come to be tied to how states treat their women, where gender quotas have been ``enacted by a range of states as a sign of democracy and a method for reducing government corruption."\footnote{  Marysia Zaleweski and Anne Runyan, “Taking Feminist Violence Seriously in International Relations,” International Feminist Journal of Politics 15(3) (2013): 293-313, p.299.}  Richard Ashley suggests that democracy can also be conceptually linked with its goals, using the democratic peace as an example.\footnote{Richard K. Ashley, “Untying the Sovereign State: A Double Reading of the Anarchy Problematique,” Millennium: Journal of International Studies 17(2) (1988): 227-262.}  Ashley explains that the ``academically certified version of the democratic peace has led to a securitization of democracy" which is deeply problematic.\footnote{Ibid.}   Andrew Linklater, on the other hand, suggests that it is the tie to Western liberalism that can be most insidious for the concept of democracy, and advocates for theorizing democracy ``without assuming that Western liberal democracy is the model of government which should apply universally."\footnote{Andrew Linklater, “The Achievements of Critical Theory,” in Ken Booth, Steve Smith, and Marysia Zalewski, eds., International Relations Theory: Positivism and Beyond (Cambridge: Cambridge University Press, 1996), 279-300, p.294.}

While many critical theorists agree that, in a variety of ways, the notion of democracy that is deployed in contemporary global politics and in contemporary IR research is both empirically and normatively problematic, they disagree strongly on how to handle it. Some suggest that the concept of democracy is now itself part of the problem (perhaps what Baudrillard in {\em The Mirror of Production} would call a repressive simulation), while others are interested in reviving a different understanding of democracy.\footnote{Jean Baudrillard, The Mirror of Production (Candor, NY: Telos Press Publishing, 1975).}   For example, Bieler and Morton suggest that the problem is not the concept of democracy itself, but the hollowing of that concept.\footnote{Andreas Bieler and Adam David Morton, “A Critical Theory Route to Hegemony, World Order, and Historical Change: Neo-Gramscian Perspectives in International Relations,” Capital \& Class 82 (2004): 85-113, p.97.}   They argue that there is a politics of supremacy that has come to replace democracy which ``involves a hollowing out of democracy and the affirmation, in matters of political economy, of a set of macro-economic policies such as market efficiency, discipline and confidence, policy credibility and competitiveness."\footnote{Ibid.}  Some critical IR theorists, then, look to rescue the concept of democracy from that hollowness.  

For example, Ken Booth (2007: 55, citing Murphy 2001: 67) suggests that what unifies critical theory in IR ``in addition to its post-Marxist sensibility, is democracy. Craig Murphy got it exactly right when he saw this emerging critical theory project being `today's manifestation of a long-standing democratic impulse in the academic study of international affairs.' In other words, it was academe’s contribution to `egalitarian practice.'"\footnote{Ken Booth, Theory of World Security (Cambridge: Cambridge University Press, 2007), p.55, citing Craig N. Murphy, “Critical Theory and the Democratic Impulse: Understanding a Century-Old Tradition,” in Richard Wyn-Jones, Critical Theory and World Politics (Boulder, CO: Lynne Rienner, 2001), pp. 61-76.}  To follow up, Booth goes over a number of different types of possible democracy, arguing that finding a good notion of democracy is key to the emancipatory mission he attributes to critical IR. Booth explains that ``there will be no emancipatory community without dialogue, no dialogue without democracy."\footnote{Ibid, p. 272.} In emancipatory critical theorists’ terms, though,  this is a different type of democracy, one that ``begins with greater recognition, representation, and access within existing institutions and demands new mechanisms for popular control of local, global, and security issues."\footnote{Alison Brysk, From Tribal Village to Global Village: Indian Rights and International Relations in Latin America (Stanford, CA: Stanford University Press, 2000), p.293.}   Following William Connolly, Richard Shapcott describes the democracy favored by critical IR as ``a democratic ethos" which is ``an ethos of pluralisation" focused on creating room for difference.\footnote{Richard Shapcott, Justice, Community, and Dialogue in International Relations (Cambridge: Cambridge University Pressm 2001), p.70.}   This leads Shapcott to express interest in ``an attempt to provide an account of democracy that does not privilege the `abstract' other and a universal subjectivity or the territorial restrictions of the nation-state."\footnote{Ibid, p. 69.}

While the dividing line is not perfect, it might be worth thinking about these differences in terms of poststructuralist and emancipatory critical theory. Both argue that there are problematic significations of current deployments of the notion of democracy in global politics. The latter is interested in reviving a more just concept of democracy, where the former is more interested in mapping the injustice that may well be inherent in the utterance and reification of the concept. What both share, in addition to critiquing current instantiations of the concept, is an interest in how democracy is being constituted, read, reproduced, and reified as a concept, both among states in global politics and among scholars of global politics interested in understanding state (and non-state) interaction. It is possible, then, to find a number of how-possible questions in critical IR analyses of democracy.  What are the conditions of possibility of current understandings of what constitutes democracy? How are the indicators of democracy that are recognized by various scholars and policy makers chosen to the exclusion of those which are not recognized? What are the relationships between various (recognized and unrecognized) indicators of democracy? How is the concept of democracy deployed (and deployable) for (and against) certain political interests? What if anything is it about the idea of democracy that allows for hollowing, encroachment, supremacy, and/or Western dominance, if such moves happen? How are relationships between the concept of democracy and its antagonists, its component parts, and/or its results formed and cemented? What is possible (or impossible) with particular conceptions of what democracy is that could change with the change (or even elimination) of the idea? 

Certainly, answers to these how-possible questions cannot be supplied easily with extant research, much less in the scope of this chapter. That said, what the how-possible questions listed above share is an interest in {\em how democracy is being read} across a variety of audiences in a variety of different ways. The remainder of this chapter suggests the plausibility and particular advantages of the formalization of concept relationships for gaining leverage on different questions about how democracy is being read. 
  
\section*{Mapping Interpretations of Democracy}
We collected more than 100 indicators used to measure democracy over eight datasets in order to gain interpretive leverage over what political scientists tend to think democracy is, how they tend to measure it, and how countries come to be classified as democracies and non-democracies.\footnote{The datasets that we compiled include the Comparative Study of Electoral Systems (\url{http://www.cses.org}, accessed 2 April 2015), Coppedge’s Democracy Diffusion data (used in Daniel Brinks and Michael Coppedge, “Diffusion is No Illusion: Neighbor Emulation in the Third Wave of Democracy,” Comparative Political Studies 39(4) (2006): 463-489); Democracy Barometer (Democracybarometer.org, accessed 2 April 2015; Wolfgang Merkel and Daniel Bochsler,  (project leaders); Karima Bousbah; Marc Bühlmann; Heiko Giebler; Miriam Hahni; Lea Heyne; Lisa Müller; Saskia Ruth; Bernhard Wessels,. Democracy Barometer. Codebook. Version 4.1. (Aarau: Zentrum für Demokratie, 2014); Varieties of Democracy (\url{https://v-dem.net}, accessed 2 April 2015; Michael Coppedge, John Gerring, Staffan I. Lindberg, Jan Teorell, David Altman, Michael Bernhard, M. Steven Fish, Adam Glynn, Allen Hicken, Carl Henrik Knutsen, Matthew Kroenig, Kelly McMann, Daniel Pemstein, Megan Reif, Svend-Erik Skaaning, Jeffrey Staton, Eitan Tzelgov, Yi-ting Wang. 2015. “Varieties of Democracy Codebook v4”. Varieties of Democracy Project: Project Documentation Paper Series; Freedom House (\url{https://freedomhouse.org/report-types/freedom-world#.VSFsXFxYyBU}, accessed 2 April 2015); Miller-Boix-Rosato Dichotomous Codings of Democracy (Carles Boix, Michael K. Miller, and Sebastian Rosato. “A Complete Data Set of Political Regimes, 1800-2007.” Comparative Political Studies 46(12): 1523-54); Pippa Norris’ Democracy Time-Series Dataset (\url{http://www.hks.harvard.edu/fs/pnorris/Data/Democracy\%20TimeSeries\%20Data/Codebook\%20for\%20Democracy\%20Time-Series\%20Dataset\%20January\%202009.pdf}, accessed 2 April 2015), Polity IV (\url{http://www.systemicpeace.org/polityproject.html}, accessed 2 April 2015), Vanhanen’s Polyarchy Dataset (\url{https://www.prio.org/Data/Governance/Vanhanens-index-of-democracy/}, accessed 2 April 2015), World Governance Indicators compiled by Daniel Kaufmann, Aart Kraay, and Massimo Mastruzzi (\url{http://info.worldbank.org/governance/wgi/index.aspx#home}, accessed 2 April 2015); and Unified Democracy Score (Daniel Pemstein, Stephen A. Meserve, and James Melton, “Democratic Compromise: A Latent Variable Analysis of Ten Measures of Regime Type,” Political Analysis 18[4] [2010]” 426-449).}  

For the purposes of our pilot analysis, we used 14 of the variables for a particular year to map country-data-points and look for commonalities. The variables that we used were from the Polity IV and Miller-Boix-Rosato datasets. From the Polity IV dataset, we used chief executive recruitment regulation, competitiveness of executive recruitment, openness of executive recruitment, executive constraints, participation regulation, participation competitiveness, as well as concept indicators for executive recruitment, executive constraint, and political competition. We used these executive-specific, individual level variables next to the Miller-Boix-Rosato macro-level variables about democratic status and change over time, including a dichotomous measure of democracy, a measure of sovereignty, a measure of democratic transition, the previous number of democratic breakdowns, and the duration of democracy in the state (consecutive years of a particular regime type). 

Geometrically, there are a number of countries that represent the same point – that is – that their values on all 14 included indicators are the same.\footnote{Some of these indicators are dichotomous, and some are on a one-to-ten scale.}   For the purposes of differing interpretations of what democracy is and the indicators of democracy, then, the countries that represent the same data point are flat:  they represent the same configuration of indicators, definitions, variables, and conclusions for the purposes of understanding the dimensionality of the concept of democracy. Using the 14 indicators that we selected in the test-year of 2007, we find 88 unique data points – that is, 88 different configurations of the 14 variables. We then looked to analyze the relationships between those data points geometrically. 

\begin{figure}
\centerline{\includegraphics[width=4in]{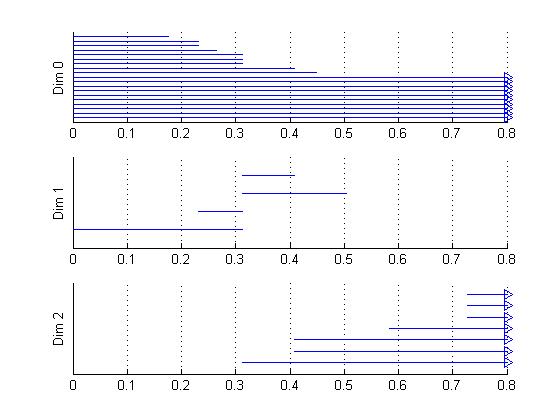}}
\caption{\label{14dbarcodes} The barcodes for the $14$-dimensional data}
\end{figure}

Using the javaplex package in Matlab, we computed the persistent homology of the Rips complex for the 88 data points. The topology of that Rips complex shows a number of distinct features. First, there are 10 connected components, represented by the U.S., Cuba, the Dominican Republic, Equatorial Guinea, Swaziland, Ivory Coast, Mauritania, Togo, Tanzania, and Guinea. A connected component is a maximal subset of a space that cannot be covered by the union of two disjoint open sets – in other words, it is a distinct and distinguishable group.  Eight of those connected components are contractible – that means that they have stronger relationships than non-contractible connected components would.\footnote{Formally, because the identity map between them is null-homotopic. As a result, in theory, a contractible space could be shrunk to a point.}  The barcodes for the relationships between datapoints can be seen in Figure \ref{14dbarcodes}. 

What these barcodes represent is the durability of certain relationships over the addition of data-points in the construction of a multi-dimensional space descriptive of the concept of democracy {\em through} collected data about its indicators.  Within the topology of the data what we analyzed, there are seven two-dimensional spheres. Six of them lie in the connected component represented by the Dominican Republic while the seventh lies in the component represented by Swaziland. The representatives of components are the states that have the most in common with the most states within the geometric component, so they matter a little as signifiers. Here, Dominican Republic is generally considered a durable if imperfect democracy, and Swaziland is generally considered autocratic.

Some of the shapes formed, then, have both clear relationships and clear implications about what ideas of democracy might be within the limitations of our truncated pilot study. Tetrahedrons formed by the Dominican Republic, Albania, Latvia, and Botswana and by the Dominican Republic, Latvia, Botswana, and Comoros are groupings of democracies with some durability and sustainability, but with (perceived) weakness on one or more indicators – here, largely, either indicators that have to do with the availability and health of political competition, or having had a democratic breakdown in the past. In this sense, they are differentiable from the democracies that constitute a single point both in theory empirically (assuming the measurements are accurate) and conceptually (that is, that there is a substantive difference about what sort of democracy a state is based on a difference on those indicators). Another tetrahedron is composed of Paraguay, Ukraine, Malawi, and East Timor. Using their Polity IV Scores (not included in the geometric analysis), those countries are considered democracies, but do not score as strongly as other countries on the democracy scale. Polity IV scores range from $-10$ (purely autocratic) to $+10$ (understood as a full democracy), but many studies that analyze democratic behavior require a score of $+6$ or higher to consider a country as a democracy. These four states often rank right around a 6 (either a little above or a little below).  These states differ from those ranked higher on the Polity IV scale primarily on three of the 14 indicators on used in this pilot study: competitiveness of participation in governance, regulation of participation in governance, and breakdowns in democratic governance.  On the other indicators, their scores are the same or substantially the same as `full' democracies. Thinking about what the components of democracy are, this suggests that democratic breakdowns distinguish countries that have had them from countries that have not, even among democracies. It also suggests that different levels of struggle with competitiveness in participation have different significations for problems with democracy.  A fourth tetrahedron is formed by Russia, Congo Kinshasa, Mozambique, and Namibia.  This tetrahedron has similar Polity IV scores to the group that was just discussed, but fares less well on a number of the indicators that we randomly selected. In the other data from the Polity Dataset, there must be a counterbalance to these countries’ negative scores not only on competitiveness and democratic breakdowns but also on openness and regulation of competition. 

\begin{figure}
\centerline{\includegraphics[width=5.5in]{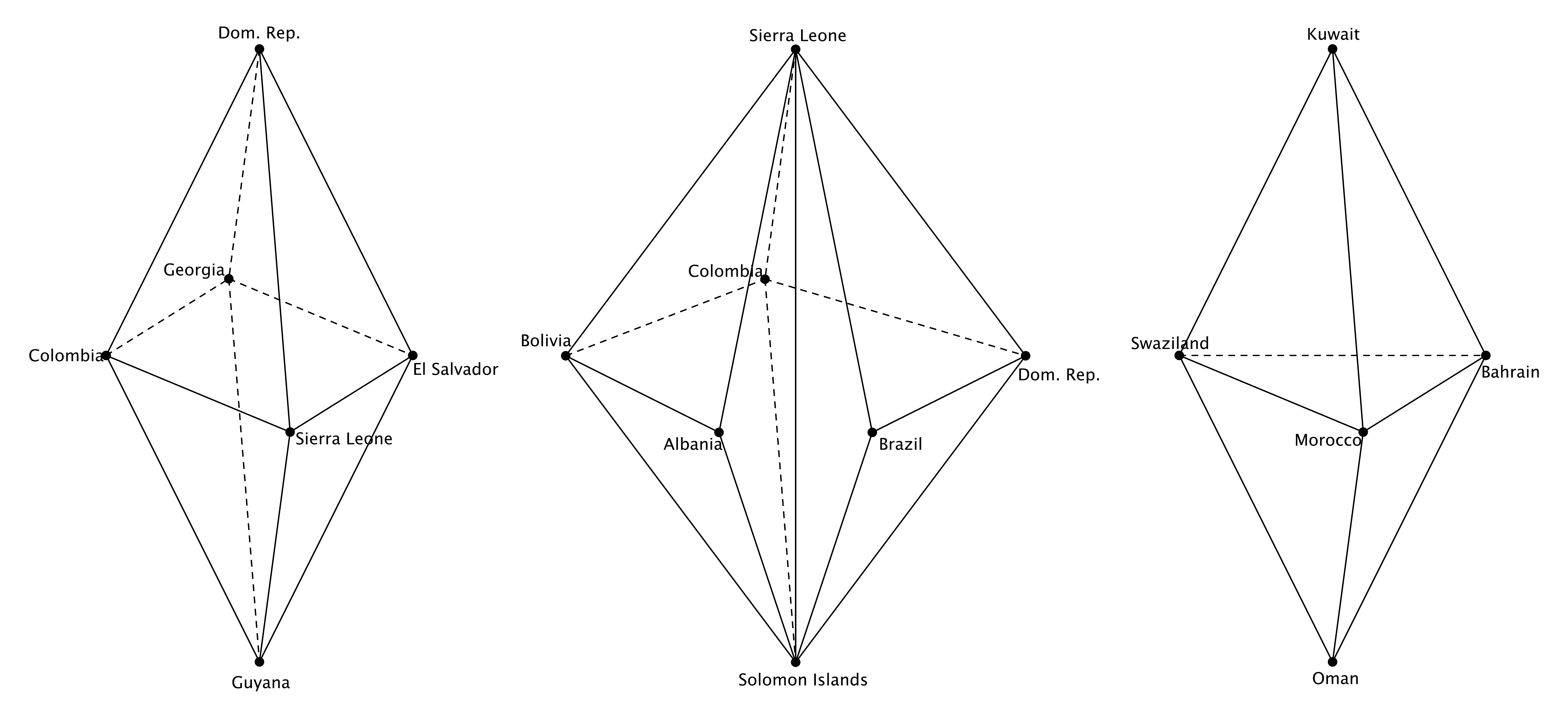}}
\caption{\label{2cycles} Some $2$-cycles in the $14$-dimensional data.  Note that these are representations of the figures in $3$-dimensional space; the actual cycles are embedded in $\zr^{14}$.}
\end{figure}

Three more complex shapes also emerge, and are depicted in Figure \ref{2cycles}.  The first is an octahedron formed by the Dominican Republic, El Salvador, Colombia, Guyana, Georgia, and Sierra Leone.  This is the figure to the left in Figure \ref{2cycles}.  These countries match each other perfectly in regulation of participation in executive elections, competitiveness of executive elections, openness of executive elections, and the competitiveness of executive recruitment. They exhibit small variations on competitiveness and regulation of participation in elections. These traits make them closely related but not collapsible into one data point. The six states differ on the existence of democratic breakdowns in their recent history, which is one of the factors that creates space among them.  What distinguishes this group as a group from other democracies is imperfect scores on participation competitiveness and regulation of executive recruitment – so it is a group generally understood to be democracies with particular weaknesses vis a vis certain indicators of democracy. 

The second more complex shape is an unnamed irregular polyhedron with ten triangular faces, composed on Dominican Republic, Colombia, Bolivia, Albania, Brazil , Solomon Islands, and Sierra Leone. This is the figure in the center of Figure \ref{2cycles}. It shares a triangular face with the octahedron on the left. This means that the two are related, but not the same.  This shape includes three of the countries in the octahedron above, with four different ones. While, in the octahedron, those countries (which score a 6 on both executive constraints indicators, that is, slightly more constrained than not) are paired with countries that score a 5 (that is, in the middle of the scale), in this irregular polyhedron, they are grouped with countries that score a 7 on those same indicators.  The combination of these shape-relationships suggests that there is both a middle ground and a threshold level for executive constraints that may be meaningful in the constitution what makes a democracy. 
	
The third more complex shape is unrelated to the first two. It is a triangular bipyramid of Swaziland, Morocco, Kuwait, Bahrain, and Oman. This is the figure on the right of Figure \ref{2cycles}.  These states rank low on most of the indicators of democracy that were included in our pilot study.  That said, their low scores vary – and are substantially less low in the area of regulation of participation in executive selection. In other words, these are countries considered autocratic that lean more in the direction of democracy when it comes to the clarity of rules of executive selection. While that does not stop them from being classified as non-democracies, it does place them in a group among those non-democracies distinguishable either from states with mediocre scores on all of the indicators or those with low scores on all of the indicators.

This is, of course, a very limited exploration of a few relationships between a few states on a few indicators over the course of one year.  And one might ask the advantage of this sort of analysis over just looking at the state’s Polity IV score, or some other aggregation of these individual indicators, rather than looking at indicators that are used to make composite measures. After all, are composite measures not their composers’ full definitions of democracy, and the (sometimes weighted) component parts their understanding of the dimensionality? Would more information not be gained by comparing composite measures, then? Certainly, such a comparison would be fruitful, and is in our future plans.\footnote{We were able to use a very limited portion of this dataset here in part because of incompatibilities across software about how to deal with missing data, and in part because of the labor-intensive nature of the initial work plan. Currently, research assistants are working on cleaning up the dataset for future work across other indicators, composite measures, and over time.}   The information such analysis would provide would be different, and no doubt an addition, as would being able to use the other hundred indicators that we collected.  The nature of this pilot study and therefore the information it can provide is very limited.  It was able to show some contours of groups of states on the basis of some commonality about particular indicators of democracy that may serve as grouping, or even tipping, points. To know {\em a meaning}, much less {\em the meaning} of democracy to IR scholars would require significantly more in-depth analysis. 

Yet such a project is well within the methodological capacity of this approach.  For example, if composite measures are by definition a combination of indicators {\em flattened} into one number, this sort of analysis can replicate those composite measures in geometric analysis so that a composite score means more than its statistical significance. It would then be possible to group states based not only on their composite scores but on the indicators on which they have the most similarity, or even on indicators on which states with widely different composite numbers share common ground. Bending those models over time can show complex configurations that may be indicators of the change of form of a state. In other words, this sort of complex mapping is capable of expanding the possibilities for categorization of states among the multiple meanings of `democracy,' and therefore of providing insight into multiple ways in which the concept is (and could be) thought about that are not always explicit in quantitative or qualitative analysis in comparative politics and IR. 

The biggest possible payoff in terms of looking at understandings of democracy, however, is the next four possible steps. First, looking at the data on democracy as abstract simplical complexes and examining the filtrations and barcodes that comes out of those constructions for the relationships {\em among states on indicators} can provide a basis for multidimensional mapping {\em of indicators in relation to each other}. In other words, it would be possible to go beyond the positivist tendency to test collinearity in order to look at the relationships among particular variables held to be indicators of democracy both definitionally and as operationalized and applied in the field. This could be done with existing data on democracy in the field of political science. A second step would be to compare the homological analysis of the indicators, measures, and definitions of democracy in political science to those used in the policy world – using critical discourse analysis\footnote{Ruth Wodak and Michael Meyer, Methods for Critical Discourse Analysis (London: Sage, 2009).}  to collect data from state policy statements, press releases, and leader quotes in news publications. A third step could be comparing that to survey data collected on the ground in a number of states around the world, where respondents were asked about what democracy is, how you can tell that a democracy exists, and whether neighboring states are democracies or not. A fourth step could compare these concept maps to concept maps for other, related but distinguishable, concepts, like those discussed by Doty, Zalewski and Runyan, and Ashley.\footnote{See notes 42-45.}  

It is possible to be critical of these ideas for topological data analysis of democracy in IR by asking – what is going on here? Is it not just a hyperactive process to critique (or even perfect) the positivist operationalization of democracy? Certainly, these methods could be used to such an end, making space for hypertechnical representations in regressions looking to figure out how regime type influences foreign policy, or in predictions about the evolutions of forms of state government. But that is neither our intent, nor, in our view, the primary benefit of this methodological innovation. 

Instead, we see the primary utility of an expansion of this sort of analysis in Shapcott’s understanding that democracy ``must forever be questioning itself and the boundaries that it invokes."\footnote{Shapcott, Justice, Community, and Dialogue, p.68.}   Mapping meanings of democracy, measurements of democracy, and comparisons on indicators of democracy in multidimensional space can help us understand the ways that various concepts are leveraged in favor of, and related to, certain notions of democracy, as well as maps and relationships of inclusion and exclusion. If it is possible to have some understanding of axes of rotation and points of engagement by just looking at a few datapoints on a few indicators in one year, the analytical possibility of a full exploration mapping understandings of democracy is almost unlimited. Such mappings could contribute not only to the analysis of some of the specific {\em how-possible} questions above, but also to other questions yet unasked interested in relationships (tensions, similarities, and the simultaneous presence of both) between different ways democracy is read in global politics and in disciplinary IR. 

\section*{The Potential Payoffs of Topological Analysis for Critical IR}
If the above formalizations of concept vectors, concept spaces, concept topologies, and contexts are applicable to any concept in any theoretical context (in the informal sense), why deploy them for use in poststructuralist IR? In other words, even if these methods {\em could work} for poststructuralist analysis, and provide some value added, why use them? Is the value-added enough? After all, even if the ontological and epistemological positions of mathematical formalism and poststructuralist IR have commonalities, those commonalities do not dictate the fruitfulness of the two working together. While those commonalities (which we point out above) are the basis for our claim that the two are compatible despite a general association of quantitative work (of whatever flavor) and positivism, and poststructuralism with qualitative methods, our argument that this analysis is usefully employed in poststructuralist analysis is more based on the capacities that these sorts of representations have that the tools traditionally available to poststructuralist scholars do not easily replace and/or replicate.

Particularly, we argue that there are three principal potential benefits to the deployment of this methodology for poststructuralist IR. The first is that the complex concept modeling can be used to reveal dimensions of concepts (formally and informally) previously underexamined, either in mainstream or in critical analysis. Thinking about the contours of concepts helps to understand not only the ideas that go into them and/or their underlying assumptions. If topological concept mapping can capture the complexity of relations between features,\footnote{See discussion in Aisbett and Gibbon, “A General Formulation of Conceptual Spaces,” p.217.}  then this provides a different way to think about the underlying assumptions, building blocks, and inscriptions, and fixings of meanings in poststructuralist terms. Multidimensional concept modeling also provides a tool to think about the change of concepts over time, over place, and in the ways that they are thought about – either in the discipline or in the policy world.

Second, and perhaps more interestingly, the tool of topological concept mapping can be based in empirical and representative studies but is not confined to them. That is, a model built to represent the dimensionality of a concept can be studied with changes to that dimensionality to see about potential changes in the concept. This serves the purposes of emancipatory critical theorizing – would the world be a better place if we thought about things differently? It also serves the purposes of post‐structuralist critical theorizing – how do concept structures become sticky and reified? What would it look like to un-stick a particular dimension of a concept? If concept mapping can be manipulated temporally,\footnote{See discussion in ibid., p.218.} utilized to analyze transitional effects,\footnote{See discussion in Ibid., p.220.}  tessellated to unpack the relations between different compiled meanings,\footnote{See discussion in Rickard, “A Conceptual Geometry for Conceptual Space,” p.315.}   morphologized to explore metaphorical relationships,\footnote{Gardenfors, “Mental Representation,” p.40.}  and translated to fuzzy geometry to understand liminality,\footnote{Rickard, “A Conceptual Geometry for Conceptual Space,” p.311.}  there is significant potential for developing critical analysis of what global politics is, how it is possible, and how it is constituted, reified, and performed.

This could be done in a way that emphasizes relative relationality – which, in our view, is the third major potential payoff for critical theorizing. While the descriptors for conceptual relationships are limited in terms of the sorts of relationships we can think about – in between, close to, far from, etc; the topological descriptors are in theory both unlimited, and more clearly specifiable (given the potential for multidimensionality). In fact, the complexity both of representation and exploration is in theory unlimited. In practice, it is limited only by the possible accessibility of time and information, and the possible specification of data points. Neither of these limits is concerning, though, given that even the least complex representations have potential exploratory value if not empirical value.

\end{document}